\documentclass[a4paper]{amsart}

\usepackage[english]{babel}

\usepackage[latin1]{inputenc}

\usepackage{amssymb}
\usepackage{amsmath}
\usepackage{amsthm}

\begin{document}

\parindent0pt

% abbreviations for \mathbb-commands:
\newcommand{\bbR}{\mathbb{R}}
\newcommand{\bbN}{\mathbb{N}}
\newcommand{\bbD}{\mathbb{D}}

% define environments:
\theoremstyle{definition}
\newtheorem{env_def}{Definition}[section]
\newtheorem{env_rem}[env_def]{Remark}
\newtheorem{env_examp}[env_def]{Example}

\theoremstyle{plain}
\newtheorem{env_prop}[env_def]{Proposition}
\newtheorem{env_lem}[env_def]{Lemma}
\newtheorem{env_theo}[env_def]{Theorem}
\newtheorem{env_cor}[env_def]{Corollary}

% title:
\title[Square Roots in Strictly Linearly Ordered Semigroups]{Square Roots and Continuity in Strictly Linearly Ordered Semigroups on Real Intervals}
\author{Jochen Gl{\"u}ck}
\date{February 27, 2014}
\address{Jochen Gl{\"u}ck, Institute of Applied Analysis, Ulm University, 89069 Ulm, Germany}
\email{jochen.glueck@uni-ulm.de}

\begin{abstract}
	In this article we show that the semigroup operation of a strictly linearly ordered semigroup on a real interval is automatically continuous if each element of the semigroup admits a square root. Hence, by a result of Acz{\'e}l, such a semigroup is isomorphic to an additive subsemigroup of the real numbers.
\end{abstract}

\maketitle

% content:

\section{Introduction and statement of the main result}

A \emph{strictly linearly ordered semigroup} is a tuple $(X, \le, \circ)$, where $X$ is a non-empty set, $\le$ is a linear order on $X$ and $\circ: X^2 \to X$ is an associative map such that $a < b$ implies $a\circ x < b \circ x$ and $x \circ a < x \circ b$ for all $a,b,x \in X$ (by $a < b$ we mean that $a \le b$, but $a \not= b$). \par
There is a long tradition of investigating the question whether a strictly linearly ordered semigroup can be embedded into the additive semigroup of all real numbers. A first result of this type was obtained by Abel in case that the semigroup operation fulfils some differentiability conditions (\cite{Abel1826}; see also \cite{Lawson1996} for a modern review of Abel's article). In \cite{Aczel1948} Acz{\'e}l gave the following generalization of Abel's result:

\begin{env_theo}[Acz{\'e}l] \label{theo_aczel}
	Let $I \subset \bbR$ be a real interval endowed with the usual order $\le$. Furthermore, let $\circ$ be an associative operation on $I$ such that $(I, \le, \circ)$ is a strictly linearly ordered semigroup. \par
	If the semigroup operation $\circ$ is jointly continuous, then there is an interval $J \subset \bbR$ and a monotone, bijective and continuous mapping $f: I \to J$ such that $f(x\circ y) = f(x) + f(y)$ for all $x,y \in I$.
\end{env_theo}

An \emph{isomorphism} of two strictly linearly ordered semigroups $(X, \le, \circ)$ and $(Y, \le, \diamond)$ is a function $f: X \to Y$ which is monotone, bijective and which fulfils the equation $f(x_1\circ x_2) = f(x_1) \diamond f(x_2)$ for all $x_1,x_2\in X$. Note that if $f: X \to Y$ is an isomorphism, then so is the inverse mapping $f^{-1}: Y \to X$. We say that $(X, \le, \circ)$ and $(Y, \le, \diamond)$ are isomorphic if there exists an isomorphism between them. \par
In this terminology Acz{\'e}l's Theorem states that if $I$ is a real interval and $(I, \le, \circ)$ is a strictly linearly ordered semigroup with jointly continuous operation $\circ$, then $(I, \le, \circ)$ is isomorphic to $(J, \le, +)$, where $J$ is a real interval, too (note that in this case the continuity of $f$, which is asserted in Acz{\'e}l's Theorem, follows automatically since every monotone bijection between two real intervals is continuous). \par

There are several generalizations of Acz{\'e}l's Theorem weakening the monotonicity condition on the binary operation (see for example \cite{Bacchelli1986} and \cite{Ling1965}) or extending the result to $n$-ary operations (see the recent article \cite{Couceiro2012}). In \cite{Craigen1989} a new, streamlined proof of Acz{\'e}l's Theorem was given by Craigen and Pal{\'e}s and in \cite{Aczel2004} Acz{\'e}l gives a review of his own proof and some notes on Craigen and Pal{\'e}s' version of the proof. The problem of embedding a linearly ordered semigroup into the real numbers under topological conditions was also treated in two recent articles by Ghiselli Ricci (\cite{GhiselliRicci2006} and \cite{GhiselliRicci2009}) \par

There are also other approaches to the subject, that use algebraical instead of topological conditions on the semigroup (see \cite{Fuchs1963}, chapter XI for an overview over some results of this type). \medskip \par

For a detailed discussion of the history of the subject see the survey article of Hofmann and Lawson in \cite{Hofmann1996}. For the connection of the subject to Hilbert's Fifth Problem, see also \cite{Aczel1989} and \cite{Hofmann1994}. \medskip \par

In this article we present a further approach to the subject discussed above, that is based on the notion of square roots: An element $x$ of a strictly linearly ordered semigroup is said to \emph{admit a square root} if there is an element $y$ such that $y^2 = y \circ y = x$. In this case $y$ is uniquely determined. It is called the \emph{square root} of $x$ and denoted by $x^{\frac{1}{2}}$. Given a strictly linearly ordered semigroup on a real interval and assuming that each element of the semigroup admits a square root, we will prove that the semigroup operation is jointly continuous, so that we can apply Acz{\'e}l's Theorem in order to embed our semigroup into the real numbers. \par
Note, that square roots (and hence, dyadic numbers) play a role in several constructions and proofs within the theory of ordered algebraic objects. See e.g. \cite{Aczel1948a}, \cite{Jager2005} and section XI.8 in \cite{Fuchs1963} for the usage of roots and dyadic numbers in some approaches to abstract mean values. Moreover, roots are involved in many of the above mentioned embedding results in order to construct some kind of exponential-functions which serve as a semigroup-homomorphisms. \par
However, while those approaches impose different topological or algebraical conditions on the semigroup and then use roots mainly for the construction of homomorphisms, we take the existence of square roots as a starting point and then derive continuity of the semigroup operation. \medskip \par

Throughout, let $I$ be a real interval, $\le$ the restriction of the usual order on $\bbR$ to $I$ and $\circ$ a binary operation on $I$ such that $(I, \le, \circ)$ is a strictly linearly ordered semigroup. \smallskip \par

Our aim is to prove the following result:
\begin{env_theo}[main result] \label{theo_main_result}
	Suppose that every element of $I$ admits a square root. Then the operation $\circ: I^2 \to I$ is jointly continuous.
\end{env_theo}

Therefore, by Acz{\'e}l's Theorem, we obtain the following corollary:
\begin{env_cor} \label{cor_isomorphism_to_the_real_line}
	Suppose that every element of $I$ admits a square root. Then $(I, \le, \circ)$ is isomorphic to the strictly linearly ordered semigroup $(J, \le, +)$, where $J$ is one of the real intervals $\bbR$, $[0, \infty)$, $(0, \infty)$ or $[0]$.
	\begin{proof}
		By Acz{\'e}l's Theorem and Theorem \ref{theo_main_result}, $(I, \le, \circ)$ is isomorphic to $(J, \le, +)$, where $J$ is a real interval. To determine the interval $J$, observe that, due to our square root condition, for each $x \in J$ the sequence $(\frac{x}{2^n})$ is contained in $J$ as well. Therefore, $0$ is in the topological closure of $J$ and since $J$ is closed under the real addition $+$, we conclude $J \in \{ \bbR, [0, \infty), (0, \infty), (-\infty,0], (-\infty,0), [0] \}$. \par
		Since the map $x \mapsto -x$ is an isomorphism between $((-\infty,0],\le,+)$ and $([0,\infty),\le, +)$ respectively between $((-\infty,0), \le, +)$ and $((0,\infty), \le, +)$, the corollary is proved.
	\end{proof}
\end{env_cor}

\section{Proof of the main result}

To prove Theorem \ref{theo_main_result}, we first recall some notions about strictly linearly ordered semigroups. \par
An element $a \in I$ is \emph{strictly positive}, if for all $x \in I$ we have $a \circ x > x$ and $x \circ a > x$. It is \emph{strictly negative}, if for all $x \in I$ we have $a \circ x < x$ and $x \circ a < x$. If for all $x \in I$ the equations $a \circ x = x \circ a = x$ hold, then $a$ is called a \emph{unit}. There is at most one unit in $I$ (for if $a$ and $\tilde a$ are units, then $a = a \circ \tilde a = \tilde a$).  Moreover, due to the following proposition, every element $a \in I$ that is not a unit, is either strictly positive or strictly negative.

\begin{env_prop} \label{prop_positive_elements}
	Let $a \in I$.
	\begin{itemize}
		\item[(i)] If $a^2 > a$, then $a$ is strictly positive. \par 
		\item[(ii)] If $a^2 < a$, then $a$ is strictly negative. \par
		\item[(iii)] If $a^2 = a$, then $a$ is a unit.
	\end{itemize}
	Suppose furthermore that there is a unit $e \in I$.
	\begin{itemize}
		\item[(iv)] The element $a$ is strictly positive if and only if $a > e$. \par
		\item[(v)] The element $a$ is strictly negative if and only if $a < e$.
	\end{itemize}
	\begin{proof}
	If $a^2 > a$ and $x \in I$, one has $a \circ (a \circ x) = a^2 \circ x > a \circ x$ and therefore $a \circ x > x$. Similarly, from $x \circ a^2 > x \circ a$ we obtain $x \circ a > x$. This shows (i). \par
	In order to show (ii) and (iii), just imitate the prove of (i), while replacing the $>$-relation by $<$ or $=$. \smallskip \par
	Now suppose that $e \in I$ is a unit. If $a \in I$ is strictly positive, then $a = a \circ e > e$. On the other hand, if we assume $a > e$ than $a \circ x > e \circ x = x$ and $x \circ a > x \circ e = x$ for each $x \in I$, i.e. $a$ is strictly positive, which proves (iv). Similarly, one proves (v).
	\end{proof}
\end{env_prop}

As a corollary we obtain:

\begin{env_cor} \label{cor_positive_elements_and_powers}
	Let $a \in I$. Then, for any $m \in \bbN$, $a^m$ is strictly positive (strictly negative, a unit) if and only if $a$ is strictly positive (strictly negative, a unit).
	\begin{proof}
		If $a$ is strictly positive, then $(a^{m})^{2} = a^{2m} > a^m$ by the definition of strict positivity. By Proposition \ref{prop_positive_elements} (i) it follows that $a^m$ is strictly positive. \par
	Similarly we can see that if $a$ is strictly negative respectively a unit, then so is $a^m$. This in turn implies that if $a^m$ is strictly positive, then $a$ can neither be strictly negative nor a unit, i.e. $a$ is strictly positive. \par
	In the same way, we see that if $a^m$ is strictly negative respectively a unit, then so is $a$.
	\end{proof}
\end{env_cor}

From now on, we suppose that every element $x \in I$ admits a square root $x^{\frac{1}{2}}$.
Let $\bbD := \{\frac{k}{2^n}: \; k \in \bbN, \; n \in \bbN_0\}$ be the set of strictly positive dyadic numbers. For each $x \in I$ and $d = \frac{k}{2^n} \in \bbD$, we define a dyadic power by $x^d = x^{\frac{k}{2^n}} = (x^{\frac{1}{2^n}})^k$, where $x^{\frac{1}{2^n}}$ is constructed with respect to the recurrence relation $x^{\frac{1}{2^n}} = (x^{\frac{1}{2^{n-1}}})^{\frac{1}{2}}$. \par
One easily verifies that $x^d$ does not depend on the representation of $d = \frac{k}{2^n}$ and is therefore well-defined. Furthermore, the usual exponential equations hold:
\begin{align*}
	x^{d_1 + d_2} & = x^{d_1} \circ x^{d_2} \quad \text{ and } \\
	(x^{d_1})^{d_2} & = x^{d_1 d_2} \text{,}
\end{align*}
where $x \in I$ and $d_1,d_2 \in \bbD$. \medskip \par

We show some elementary properties of the powers $x^d$ in case when $x$ is strictly positive:
\begin{env_prop} \label{prop_monotone_dyadic_powers}
	Let $x \in I$ be strictly positive.
	\begin{itemize}
		\item[(i)] If $d \in \bbD$, then $x^d$ is strictly positive, too. \par
		\item[(ii)] The map $\bbD \to I$, $d \mapsto x^d$ is strictly increasing.
	\end{itemize}
	\begin{proof}
		(i) Let $d = \frac{k}{2^n} \in \bbD$. By Corollary \ref{cor_positive_elements_and_powers}, $x^k$ is strictly positive. However, if $x^{\frac{k}{2^n}}$ was strictly negative or a unit, then, again by Corollary \ref{cor_positive_elements_and_powers}, $x^k = (x^{\frac{k}{2^n}})^{2^n}$ would be strictly negative or a unit, too. Hence, $x^{\frac{k}{2^n}}$ is strictly positive. \par
		(ii) Let $d_1,d_2 \in \bbD$ and $d_1 < d_2$. It follows from (i) that $x^{d_2 - d_1}$ is strictly positive, so $x^{d_2} = x^{d_2 - d_1} \circ x^{d_1} > x^{d_1}$.
	\end{proof}
\end{env_prop}

Of course, similar results hold if $x$ is strictly negative. \medskip \par

Let $x \in I$ be strictly positive. We define a map
\begin{align*}
	\exp_x: \; \bbR_{>0} \to I \text{,} \qquad r \mapsto \inf\{x^d: \; d \in \bbD, \; d \ge r\} \text{.}
\end{align*}

Note that $\exp_x(r)$ is well-defined, for if $r > 0$, then the set $\{x^d: \; d \in \bbD, \; d \ge r\}$ is bounded below by $x^{\tilde d} \in I$ for a sufficiently small $\tilde d \in \bbD$. Thus, the infimum in the above definition exists and is contained in $I$. \par
\begin{env_prop} \label{prop_real_exponential}
	Let $x \in I$ be strictly positive. Then the map $exp_x$ is strictly increasing, and for each $d \in \bbD$ the equation $\exp_x(d) = x^d$ holds.
	\begin{proof}
		Let $d \in \bbD$. By Proposition \ref{prop_monotone_dyadic_powers} we have $x^d = \min\{x^{\tilde d}: \; \tilde d \in \bbD, \; \tilde d \ge d\}$, which implies $\exp_x(d) = x^d$. \par
		Now, let $r_1, r_2 \in \bbR_{>0}$ and $r_1 < r_2$. There are $d_1,d_2 \in \bbD$ such that $r_1 < d_1 < d_2 < r_2$. Thus, we obtain
		\begin{align*}
			\exp_x(r_1) \le x^{d_1} < x^{d_2} \le \exp_x(r_2) \text{,}
		\end{align*}
		in which the latter two inequalities follow from Proposition \ref{prop_monotone_dyadic_powers}.
	\end{proof}
\end{env_prop}

According to the above proposition, we can think of $\exp_x$ as an extension of the map $d \mapsto x^d$ to the strictly positive real line. However, since we have not yet proved continuity of $\circ$, it is not clear at all, whether the homomorphism properties of the dyadic powers $x^d$ are fulfilled by the real powers $\exp_x(r)$, too. \smallskip \par

Nevertheless the exponential map $\exp_x$ is quite useful in the proof of the following lemma, which is the essential step on our way to prove Theorem \ref{theo_main_result}:

\begin{env_lem} \label{lem_approx_unit}
	Let $x \in I$ be strictly positive. Then for each $a \in I$, we have $\inf_{n \in \bbN} a \circ x^{\frac{1}{2^n}} = \lim_{n \to \infty} a \circ x^{\frac{1}{2^n}} = a$.
	\begin{proof}
		By Proposition \ref{prop_monotone_dyadic_powers} the sequence $a \circ x^{\frac{1}{2^n}}$ is decreasing and bounded below by $a$, so we have $b := \inf_{n \in \bbN} a \circ x^{\frac{1}{2^n}} = \lim_{n \to \infty} a \circ x^{\frac{1}{2^n}} \ge a$. \par
Assume towards a contradiction that $b > a$. Then for each $r > 0$ the open interval $J_r := (a \circ \exp_x(r), \; b \circ \exp_x(r))$ is non-empty (since $I$ is a real interval and therefore densely ordered). What is more, all the intervals $J_r$ are pairwise disjoint. To see this, let $r_1,r_2$ be distinct elements of $\bbR_{>0}$, say $r_1 < r_2$. Choose $d_1,d_2 \in \bbD$ such that $r_1 < d_1 < d_2 < r_2$ and choose $m \in \bbN$ such that $\frac{1}{2^m} \le d_2 - d_1$. Then we obtain
		\begin{align*}
& a \circ \exp_x(r_2) \ge a \circ \exp_x(d_2) = a \circ x^{d_2} = a \circ x^{d_2-d_1} \circ x^{d_1} \ge \\
& \ge a \circ x^{\frac{1}{2^m}} \circ \exp_x(d_1) \ge b \circ \exp_x(d_1) \ge b \circ \exp_x(r_1) \text{.}
		\end{align*}
		Thus, $J_{r_1} \cap J_{r_2} = \emptyset$. However, this is a contradiction: Since $I$ is a real interval, it is separable (i.e. there is a countable, dense subset $A \subset I$) and thus there is no uncountable family of pairwise disjoint, non-empty open intervals on $I$. \par
		This implies that our assumption is false and therefore $b = a$.
	\end{proof}
\end{env_lem}

In the setting of a real interval, separability arguments are often useful to prove continuity results (see, for example, \cite{Ruppert1984}, section II.5.3, where a similar argument is used to deduce, under certain conditions, the joint continuity of a semigroup operation from the continuity of all left translations $l_a: x \mapsto a \circ x$). \smallskip \par

From Lemma \ref{lem_approx_unit} we obtain the following corollary:

\begin{env_cor} \label{cor_continuous_1}
	For the left translations $l_a: I \to I$, $x \mapsto a \circ x$ ($a \in X$), the following assertions hold:
	\begin{itemize}
		\item[(i)] If there is at least one strictly positive element in $I$, then all left translations $l_a$ are continuous from above on $I$. \par
		\item[(ii)] For each strictly positive element $x_0 \in X$, all left translations $l_a$ are continuous from below at $x_0$.
	\end{itemize}
	\begin{proof}
		Let $x_0 \in I$ be strictly positive, let $x \in I$ and let $(x_n) \subset I$ be a sequence converging to $x$, where $x_n > x$ for each $n \in \bbN$. If $m \in \bbN$, we have $x < x_n < x \circ x_0^{\frac{1}{2^m}}$ for sufficiently large $n$. Thus, by Lemma \ref{lem_approx_unit},
		\begin{align*}
			a \circ x < a \circ x_n < a \circ x \circ x_0^{\frac{1}{2^m}} \to a \circ x \quad (m \to \infty) \text{.}
		\end{align*}
		This shows that $l_a$ is continuous from above at each $x \in I$. \par
		Now, consider the strictly increasing sequence given by $y_n := x_0^{1 - \frac{1}{2^n}} < x_0$. We show that $a \circ y_n \to a \circ x_0$ for each $a \in I$: \par
		Let $b < a \circ x_0$. By Lemma \ref{lem_approx_unit} we have $b \circ x_0^{\frac{1}{2^{n_0}}} < a \circ x_0$ for some $n_0 \in \bbN$. Then for each $n \ge n_0$, the inequalities
		\begin{align*}
			b \circ x_0^{\frac{1}{2^{n_0}}} < a \circ x_0  \le a \circ y_n \circ x_0^{\frac{1}{2^{n_0}}}
		\end{align*}
		hold. Therefore $b < a \circ y_n < a \circ x_0$, so we showed that $a \circ y_n$ converges to $a \circ x_0$. As in the first part of this proof, this implies $a \circ x_n \to a \circ x_0$ if $x_n < x_0$ is a sequence converging to $x_0$ from below. Hence, $l_a$ is continuous from below at $x_0$.
	\end{proof}
\end{env_cor}

In order to further strengthen the continuity result of Corollary \ref{cor_continuous_1} we make use of the following observations:
\begin{env_rem} \label{rem_associated_semigroups}
	From $(I, \le, \circ)$ we construct two related strictly linearly ordered semigroups:
	\begin{itemize}
		\item[(i)] Consider the binary operation $\hat \circ: I^2 \to I$ that is defined by $a \; \hat \circ \; b := b \circ a$. Then $(I, \le, \hat \circ)$ is a strictly linearly ordered semigroup that fulfils all conditions that we have required for $(I, \le, \circ)$ at the beginning. \par
		\item[(ii)] Let $\tilde I$ be a real interval (with the usual order $\le$) such that there is a decreasing, bijective map $\varphi: I \to \tilde I$. We define a binary operation $\diamond$ on $\tilde I$ by setting $\tilde a \diamond \tilde b := \varphi(\varphi^{-1}(\tilde a) \circ \varphi^{-1}(\tilde b))$ for all $\tilde a, \tilde b \in \tilde I$. Then $(\tilde I, \le, \diamond)$ is a strictly linearly ordered semigroup that fulfils all condition that we required for $(I, \le, \circ)$ at the beginning of the paper. Furthermore, observe that $\varphi$ is a homeomorphism (since $\varphi$ and $\varphi^{-1}$ are bijective, decreasing maps between two real intervals).
	\end{itemize}
\end{env_rem}

\begin{env_cor} \label{cor_continuous_2}
	For each $a \in I$ the left translation $l_a: \; x \mapsto a \circ x$ and the right translation $r_a: \; x \mapsto x \circ a$ are continuous on $I$.
	\begin{proof}
		By part (i) of the preceding remark it suffices to consider all left translations $l_a$. For if the assertion is proved for all left translations, then all left translations in $(I, \le, \hat \circ)$ are continuous as well, and those are exactly the right translations in $(I, \le, \circ)$. \smallskip \par
		We may suppose, that $I$ contains more than one element and we distinguish between several cases in order to reduce our assertion to Corollary \ref{cor_continuous_1}:
		\begin{itemize}
			\item[(a)] If $I$ contains only strictly positive elements, then the assertion is immediate from Corollary \ref{cor_continuous_1}. Now, assume that there are only strictly negative elements in $(I, \le, \circ)$. Define another strictly linearly ordered semigroup $(\tilde I, \le, \diamond)$ with a decreasing semigroup-isomorphism $\varphi: I \to \tilde I$ as in Remark \ref{rem_associated_semigroups} (ii). Then $\tilde I$ consists only of strictly positive elements. This implies that every left translation in $\tilde I$ is continuous and since $\varphi$ is a homeomorphism, every left translation in $I$ is continuous, too. \par
			\item[(b)] If $I$ contains a unit and all other elements of $x$ are strictly positive, then each left translation $l_a$ is continuous from above on $I$ and continuous from below at each strictly positive element, due to Corollary \ref{cor_continuous_1}. However, by Proposition \ref{prop_positive_elements} (iv), there is no element that is lower then the unit. Hence, $l_a$ is continuous from below at the unit, too. As before, we may also reduce to this case, if $I$ contains a unit and all other elements are strictly negative. \par
			\item[(c)] Finally, assume that $I$ contains strictly positive elements as well as strictly negative elements. By Corollary \ref{cor_continuous_1}, each left translation $l_a$ is continuous from above on $I$. However, let $(\tilde I, \le, \diamond)$ and $\varphi: I \to \tilde  I$ be defined as in Remark \ref{rem_associated_semigroups} (ii), again. Then the left translation $l_{\varphi(a)}$ is continuous from above on $\tilde I$. Since $\varphi$ is a decreasing homeomorphism, this implies that $l_a$ is continuous from below on $I$. \hfill \qedhere
		\end{itemize}
	\end{proof}
\end{env_cor}

The preceding Corollary shows that the operation $\circ$ is separately continuous. To complete the proof of Theorem \ref{theo_main_result} we have to show that this implies joint continuity of $\circ$:

\begin{proof}[Proof of Theorem \ref{theo_main_result}]
	First, note that if $U$ is a neighbourhood of an element $x \in I$, then there is a closed interval $[a,b] \subset U$ which is also a neighbourhood of $x$. \par
	Now, let $x_1, x_2 \in I$ and let $[a,b]$ be a neighbourhood of $x_1 \circ x_2$. By $\operatorname{Int}{[a,b]}$ we denote the interior of $[a,b]$. \par
	Due to Corollary \ref{cor_continuous_2}, the right translation $r_{x_2}$ is continuous and thus, there is a neighbourhood $[a_1,b_1]$ of $x_1$ such that $r_{x_2}([a_1,b_1]) \subset \operatorname{Int}{[a, b]}$. Since the left translations $l_{a_1}$ and $l_{b_1}$ are continuous as well, there is a neighbourhood $[a_2,b_2]$ of $x_2$ such that $l_{a_1}([a_2, b_2]) \subset [a, b]$ and $l_{b_1}([a_2,b_2]) \subset [a, b]$. \par
	By a monotonicity argument, this implies $y_1 \circ y_2 \in [a, b]$ for each tuple $(y_1,y_2) \in [a_1,b_1] \times [a_2, b_2]$. This completes our proof, since $[a_1, b_1] \times [a_2, b_2]$ is a neighbourhood of $(x_1, x_2)$ in $I \times I$.
\end{proof}

The preceding proof of course shows, that separate continuity implies joint continuity in every strictly linearly ordered semigroup. This and some further results of this type can be found in \cite{Ruppert1984}, Section II.5. \par

\begin{env_rem} \label{rem_exp_is_homomorphism}
	Using the continuity of $\circ$, it is not difficult to prove that the exponential map $\exp_x$ (where $x$ is strictly positive) is a continuous semigroup homomorphism that maps $\bbR_{>0}$ onto the set of all strictly positive elements of $I$. Thus, if $I$ consists only of strictly positive elements, we may derive Corollary \ref{cor_isomorphism_to_the_real_line} from our main result without using Acz{\'e}l's Theorem. 
\end{env_rem}

\section{Final remarks}

Finally, we shall discuss our condition on $I$ to be a real interval. In order to do so, recall the following order theoretical characterization of real intervals:
\begin{env_prop} \label{prop_characterization_of_real_intervals}
	Let $X$ be a non-empty set endowed with a linear order $\le$. Then $X$ is order isomorphic to a real interval $J$ (i.e. there is a monotone, bijective map $f: X \to J$) if and only if the following conditions are fulfilled:
	\begin{itemize}
		\item[(i)] The linear order $\le$ is conditionally complete on $X$, i.e. every non-empty subset $A \subset X$ that is bounded below, has an infimum in $X$. \par
		\item[(ii)] The order $\le$ on $X$ is dense, i.e. if $x_1,x_2 \in X$ and $x_1 < x_2$, then we can find an element $y \in X$ satisfying $x_1 < y < x_2$. \par
		\item[(iii)] The order topology on $X$ is separable.
	\end{itemize}
	\begin{proof}
		Cf. Definition 2.28 and Theorem 2.30 in \cite{Rosenstein1982}.
	\end{proof}
\end{env_prop}

Proposition \ref{prop_characterization_of_real_intervals} shows that there is no order theoretical distinction between real intervals and linearly ordered sets with dense, conditionally complete order and separable order topology. \par 
Therefore, it makes no difference to replace our real interval $I$ throughout the paper by a linearly ordered set $(X,\le)$ that satisfies the three conditions (i) - (iii) from Proposition \ref{prop_characterization_of_real_intervals}. We shall now discuss, whether our results remain true, if we omit one of these three conditions. \smallskip \par

Especially the role of the separability condition is interesting: At the end of their simplified proof of Acz{\'e}l's Theorem in \cite{Craigen1989}, the authors point out that their proof, though worked out only on real intervals in their article, works as well on any other set with a dense and conditionally complete order. Therefore, one can replace the real interval $I$ in Acz{\'e}l's Theorem \ref{theo_aczel} by a linearly ordered set with dense and conditionally complete order, and the assertion remains true. \par
In contrast to this situation, separability is a crucial condition in our proof of Lemma \ref{lem_approx_unit}. Indeed, if we replace the real interval $I$ by a conditionally order complete and densely ordered set $(X,\le)$ whose order topology is not separable, then Lemma \ref{lem_approx_unit} and our main result may fail, as the following example shows:

\begin{env_examp} \label{examp_inseparability}
	Let $X$ be the disjoint union $\dot{\bigcup}_{r \in \bbR_{\ge 0}} [0,r]$. We may denote each element of this union by a tuple $(r,x)$, where $r \ge 0$ and $x \in [0,r]$. Now, let $\le$ be the lexicographical order on $X$, i.e. $(r_1,x_1) \le (r_2, x_2)$ if either $r_1 < r_2$ or if $r_1 = r_2$ and $x_1 \le x_2$. Furthermore, we endow $X$ with a semigroup operation $+$, defined by $(r_1, x_1) + (r_2, x_2) = (r_1 + r_2, x_1 + x_2)$. \par
	With this definition, one immediately checks that $(X, \le, +)$ is a strictly linearly ordered semigroup. Moreover, the following holds:
	\begin{itemize}
		\item[(i)] The order $\le$ on $X$ is dense and conditionally complete. Furthermore, each element $(r,x)$ admits square root (or, using an additive notion, a half element) $(r/2, x/2)$. Nevertheless, the order topology on $X$ is not separable. So, our strictly linearly ordered semigroup $(X,\le,+)$ satisfies all conditions required in this paper, except separability. \par
		\item[(ii)] The assertion of Lemma \ref{lem_approx_unit} does not hold for $X$. To see this, consider the sequence $(\frac{1}{2^n}, 0) \in X$ which may be written as $x^{\frac{1}{2^n}}$ with $x = (1,0)$. Then we have $(1,0) + (\frac{1}{2^n},0) \to (1,1) \not= (1,0)$ (here, convergence is understood with respect to the order topology, of course). \par
		\item[(iii)] The assertion of Theorem \ref{theo_main_result} also fails in this situation: Observe that $(\frac{1}{2^n},0) \to (0,0)$, but by (ii) we have $(1,0) + (\frac{1}{2^n},0) \to (1,1) \not= (1,0) + (0,0)$. Thus the left translation $l_{(1,0)}$ is discontinuous from above in $(0,0)$.
	\end{itemize}
\end{env_examp}

A similar construction as in Example \ref{examp_inseparability} shows that we cannot drop the condition of conditionally order completeness on our linearly ordered set (consider the disjoint union $X := \dot{\bigcup}_{d \in \bbD \cup \{0\}}[0,d]$ and make similar observations as we made in Example \ref{examp_inseparability}). \par
However, it is easy to see that the density condition for $\le$ can be omitted: If $x < y$ we have $x^2 < x \circ y < y^2$ and therefore $x < (x \circ y)^{\frac{1}{2}} < y$. Thus the order $\le$ is automatically dense. \smallskip \par

Finally, the reader might have noticed, that it is not important for the roots in our semigroup to actually be square roots. One can replace our square root condition by the condition that each element $x$ admits a $k$-th root, where $k \ge 2$ is any fixed integer that is independent of $x$.

% billiography:
\bibliographystyle {plain}
\bibliography{}

\end{document}